\documentclass[12pt]{article}
\usepackage{amsfonts}
\usepackage{amssymb}
\usepackage{amsthm}
\usepackage{amsmath}
\usepackage{latexsym}
\usepackage{color}
\begin{document}
\newtheorem{lemma}{Lemma}
\newtheorem{pron}{Proposition}
\newtheorem{re}{Remark}
\newtheorem{thm}{Theorem}
\newtheorem{Corol}{Corollary}
\newtheorem{exam}{Example}
\newtheorem{defin}{Definition}

\newcommand{\la}{\frac{1}{\lambda}}
\newcommand{\sectemul}{\arabic{section}}
\renewcommand{\theequation}{\sectemul.\arabic{equation}}
\renewcommand{\thepron}{\sectemul.\arabic{pron}}
\renewcommand{\thelemma}{\sectemul.\arabic{lemma}}
\renewcommand{\there}{\sectemul.\arabic{re}}
\renewcommand{\thethm}{\sectemul.\arabic{thm}}
\renewcommand{\theCorol}{\sectemul.\arabic{Corol}}
\renewcommand{\theexam}{\sectemul.\arabic{exam}}
\renewcommand{\thedefin}{\sectemul.\arabic{defin}}
\def\REF#1{\par\hangindent\parindent\indent\llap{#1\enspace}\ignorespaces}
\def\lo{\left}
\def\ro{\right}
\def\be{\begin{equation}}
\def\ee{\end{equation}}
\def\beq{\begin{eqnarray*}}
\def\eeq{\end{eqnarray*}}
\def\bea{\begin{eqnarray}}
\def\eea{\end{eqnarray}}
\def\d{\Delta_T}
\def\r{random walk}
\def\o{\overline}

\title{The local asymptotic estimation for the supremum of a random walk with generalized strong subexponential summands
\thanks{Research supported by the National Natural Science Foundation of China
(No.11071182 ), the National Natural Science Foundation of China
(No.11401415), the Natural Science Foundation of the Jiangsu Higher
Education Institutions of China (No.13KJB110025).}}
\author{\small
Yuebao Wang$^{1}$\thanks{Corresponding author.
Telephone: 86 512 67422726. Fax: 86 512 65112637. E-mail:
ybwang@suda.edu.cn}~~ Hui Xu$^{1}$~~  Dongya Cheng$^{1}$~~ and~~ Changjun Yu$^{2}$
\\
{\footnotesize\it 1. School of Mathematical Sciences, Soochow
University, Suzhou 215006, China}\\ {\footnotesize\it 2. School of
Sciences, Nantong University, Nantong 226019, China}}
\date{}

\maketitle {\noindent\small {\bf Abstract }}\\

{\small In this paper, the local asymptotic estimation for the supremum of a random walk and its applications are presented. The summands of the random walk have common long-tailed and generalized strong subexponential distribution. This distribution class and the corresponding generalized local subexponential distribution class are two new distribution classes with some good properties. Further, some long-tailed distributions with intuitive and concrete forms are found, which show that the intersection of the two above-mentioned distribution classes with long-tailed distribution class properly contain the strong subexponential distribution class and the locally subexponential distribution class, respectively.\\

} \\


\section{Introduction}

In this paper, we primarily study the local asymptotics
for the supremum of a random walk generated by summands with
common long-tailed and generalized strong subexponential distribution.
The generalized strong subexponential distribution class and
corresponding generalized local subexponential distribution class
are two new distribution classes with some good properties.
Therefore, we introduce some related concepts and notations in this section.

Unless otherwise stated, we always assume that a
random variable (r.v.) $X$ has a proper
distribution $F$ supported on $[0,\infty)$, that is its tail distribution function
$\overline{F}(x)=P(X>x), x\in (-\infty,\infty)$ is always positive. By definition, a
distribution $F$ is said to be heavy-tailed, if for all $\alpha>0$,
\begin{eqnarray*}
\int_{0}^{\infty}e^{\alpha y}dF(y)=\infty.
\end{eqnarray*}
Otherwise, $F$ is called light-tailed. As is known to all,
heavy-tailed distributions have important applications in various
fields of applied probability, such as risk
theory, queuing system, warehousing management, branching theory,
communication net and infinite divisible distribution theory, see,
for example, Embrechts et al. (1997) and Foss et al. (2013).
Further, some common heavy tailed distributions are also studied and applied
in the statistics, see Zeller et al. (2012), Tavangar and Hashemi (2013), Sultan and Al-Moisheer (2013),
and so on. So they attract much interest of the researchers. However, the
heavy-tailed distribution class is too large and it contains some
distributions which can not be ``dominated'', so some subclasses of
the heavy-tailed distribution class with good properties were
introduced. Here we first recall some existing subclasses of the
heavy-tailed distribution class, and then introduce some new ones.


We say that a distribution $F$ belongs to the subexponential distribution class, denoted
by $F\in{\cal{S}}$, if
$$\overline{F^{*2}}(x)\sim 2\overline{F}(x),$$
that is $\lim\overline{F^{*2}}(x)(\overline{F}(x))^{-1}=2$, here and after, all limits refer to $x$ tending to infinity, unless otherwise stated. The
subexponential distribution class was introduced by Chistyakov
(1964) in the study of the branching process, where it was proved
that the subexponential distribution class is contained in the following
heavy-tailed distribution subclass. We say that a distribution $F$ belongs to
the long-tailed distribution class, denoted by $F\in{\cal{L}}$, if
for any $y\in(-\infty,\infty)$,
$$\overline{F}(x+y)\sim\overline{F}(x). $$

For convenience, we write
\begin{eqnarray*}
C_*(F)=\liminf\overline{F^{*2}}(x)(\overline{F}(x))^{-1}{\rm{~and~}}C^*(F)
=\limsup\overline{F^{*2}}(x)(\overline{F}(x))^{-1}.
\end{eqnarray*}
It was obtained by Theorem 1 of Foss and Korshunov (2007) that for any
heavy-tailed distribution $F$, $C_*(F)=2$. And it is obvious that
$F\in{\cal{S}}$ {if and only if} $C_*(F)=C^*(F)=2$, which means that
for a subexponential distribution $F$, compared with other
distributions, the fluctuations of the ratios
$\overline{F^{*2}}(x)\big(\overline{F}(x)\big)^{-1}$ is minimal as
$x\to\infty$. So to some extent, we may regard  such distribution
$F$ as ``optimal". Also, we may say a distribution $F$ is
``controllable" if $$2\le C^*(F)<\infty.$$ And when $C^*(F)=\infty$, we
say that the distribution $F$ is ``uncontrollable".

In a probabilistic model with heavy-tailed distributions, if we may
choose distributions freely, then the first choice is of course the
subexponential distributions. However, due to the complexity of the real world, the
distributions usually are not decided by us. So it is necessary to study
``controllable" distributions, even ``uncontrollable" distributions.


Kl\"{u}ppelberg (1990) first considered the ``controllable" distributions and called them ``weak idempotents". Shimura and
Watanabe (2005) called the distributions
generalized subexponential and denoted the class of such
distributions by $\cal{OS}$. In the terminology of Bingham et al.
(1987), the distributions from the class $\cal{OS}$ are called
O-regularly varying functions. Here we continue to use the notation
$\cal{OS}$.

The distribution class $\cal{OS}$ is a rather large class {which} contains many
heavy-tailed and light-tailed distributions. For research on $\cal{OS}$,
besides the above-mentioned literature, the
reader can refer to Kl\"{u}ppelberg and Villasenor (1991), Watanabe and
Yamamura (2010), Lin and Wang (2012), Cheng and Wang (2012), Cheng
et al. (2012), Yu and Wang (2013), Gao et al. (2013),
Beck et al. (2013), Xu et al. (2015a,b,c), and so on.

In present paper, we are interested in some subclasses of the class
${\cal{OS}}$, which respectively correspond two subclasses of the class
${\cal{S}}$ as follows.

We say that a distribution $F$ belongs to the strong subexponential distribution class,
denoted by $F\in{\cal{S}}^*$, if $0<EX<\infty$ and
$$\int_{0}^{x}\overline{F}(x-y)\overline{F}(y)dy\sim 2
EX\overline{F}(x).$$
In the above formulas,
if the distribution $F$ is supported on $(-\infty,\infty)$, then
$EX$ is replaced by $EX^{+}$, where $X^{+}=X\textbf{1}(X\ge0)$.

We say that a distribution $F$ belongs to the locally long-tailed distribution class, denoted by
$F\in{\cal{L}}_{\Delta_T}$, {where $T$ is some positive constant or
$\infty$}, if for some constant $x_0>0$, $F(x+\Delta_T)=P(X\in
{x+\Delta_T})>0$ for all $x\ge x_0$ and the relationship
$$F(x+y+\Delta_T)\sim F(x+\Delta_T)$$ holds uniformly for all
$y\in(0,T]$, where $\Delta_T=(0,T],~x+\Delta_T=(x,x+T]$ when
$T<\infty$, and $\Delta_T=(0,\infty),~x+\Delta_T=(x,\infty)$
when $T=\infty$. Further, we say that a distribution $F$ belongs to the locally subexponential distribution class, denoted by $F\in{\cal{S}}_{\Delta_T}$, if for some
$0<T\le \infty$, $F\in{\cal{L}}_{\Delta_T}$ and
$$F^{*2}(x+\Delta_T)\sim 2F(x+\Delta_T).$$

The strong subexponential distribution and the
locally subexponential distribution were introduced by
Kl\"{u}ppelberg (1988) and Asmussen et al. (2003), respectively
. Inspired by the distribution classes $\cal{S}^*$, ${\cal{S}}_{\Delta_T}$ and
$\cal{OS}$, we introduce the following two new distribution
classes which are the main object of study of present paper.

\begin{defin}
We say that a distribution $F$ belongs to the generalized strong subexponential distribution class, denoted by $F\in{\cal OS}^*$, if
$$C^\otimes(F)=\limsup\int_{0}^{x}\overline{F}(x-y)\overline{F}(y)dy\big(\overline {F}(x)\big)^{-1}<\infty.$$
\end{defin}

\begin{defin}
We say that a distribution $F$ belongs to the generalized locally subexponential distribution class for some
$0<T\le \infty$, denoted by $F\in{\cal OS}_{\Delta_T}$, if
for some constant $x_0>0$, $F(x+\Delta_T)>0$ for all $x\ge x_0$ and
$$C^T(F)=\limsup F^{*2}(x+\Delta_T)\big(F(x+\Delta_T)\big)^{-1}<\infty.$$
\end{defin}

Obviously, like the class ${\cal{OS}}$, the classes ${\cal{OS}}^*$ and
${\cal{OS}}_{\Delta_T}$ contain many heavy-tailed distributions
and light-tailed distributions. Moreover, the classes ${\cal{OS}}^*$
and ${\cal{OS}}_{\Delta_T}$ have also certain ``controllability''.
Existing research and application on class $\mathcal{OS}^*$, for example, can be found in Proposition 1.1, 1.2, 1.4 and Lemma 2.1 of Xu et al. (2015b).

Similar to $C_*(F)$, we write
\begin{eqnarray*}
C_\otimes(F)=\liminf\int_0^x\overline{F}(x-y)\overline{F}(y)dy(\overline{F}(x))^{-1}
\end{eqnarray*}
and
\begin{eqnarray*}
C_T(F)=\liminf F^{2*}(x+\d)(F(x+\d))^{-1}.
\end{eqnarray*}

For a heavy-tailed distribution $F$, apart from proving the fact
that $C_*(F)=2$, Foss and Korshunov (2007) also proved  that
$C_{\otimes}(F)=2EX$. However, for a light-tailed distribution $F$,
the equalities $C_{\otimes}(F)=2EX$ do not necessarily hold,
see Foss and Korshunov (2007). Similarly, for a locally heavy-tailed distribution $F$,
the equality $C_{T}(F)=2$ for some $0<T<\infty$ also does not necessarily hold, but if $F\in{\cal{L}}_{\Delta_T}$ for some $0<T<\infty$, then $C_{T}(F)=2$, see Proposition 4.1 and Remark 4.1 of Chen et al. (2013).


The remainder of this paper consists of four sections. In Section 2,
the relationships among the two new distribution classes and some
existing related ones are discussed. Some examples of long-tailed distribution show that
the class $\mathcal{OS}^*$ and
the class $\mathcal {OS}_{\Delta_T}$ properly contain the
class $\mathcal{S}^*$ and the class $\mathcal {S}_{\Delta_T}$, respectively. It should be said that the methods of construction of these distributions are not trivial, but these distributions are not particularly weird, especially their integral tail distributions are much more normal. And proofs of this examples are given in the Section 5. In Section 3, the local asymptotic estimation for the supremum of a
random walk is presented, where the summands of the random walk have
common long-tailed and generalized strong subexponential distribution. To this end, we find out some relationship between the random walk with heavy-tailed summands and the random walk with light-tailed summands. Some applications of the above results are given in Section 4.


\section{The relationships among the distribution classes}
\setcounter{equation}{0}
\setcounter{thm}{0}\setcounter{Corol}{0}\setcounter{lemma}{0}

\subsection{ The relation between the classes
$\cal{L}\cap{\cal{OS}^*}$ and $\cal{S}^*$.}

\begin{pron}\label{p201}
The inclusion relation $\cal{S}^*\subset\cal{L}\cap{\cal{OS}^*}$ is
proper.
\end{pron}

\proof Obviously, the inclusion relation $\cal{S}^*\subset\cal{L}\cap{\cal{OS}^*}$ holds. So, we just prove that the relationship is proper through the following two types of distributions.

\begin{exam}\label{e201}
Let $m\geq 1$ be any integer. Choose any constant $\alpha\in (m^{-1},
1+m^{-1})$ and any constant $x_1>4^{m\alpha(m\alpha-1)^{-1}}$. For all
integers $n\geq1$, let $x_{n+1}=x_n^{2-(m\alpha)^{-1}}$. Clearly,
$x_{n+1}>4x_n$ and $x_n\to\infty$ as $n\to\infty$. Now, define the
distribution $F$ as follows:

\begin{eqnarray}
\overline{F}(x)&=&\textbf{\emph{1}}(x< 0)+(x_1^{-1}(x_1^{-\alpha}-1)x+1)
\textbf{\emph{1}}(0\leq x< x_1)\nonumber\\
&&\ \ \ \ \ \ \ \ \ +\sum\limits_{n=1}^{\infty}
\Big(\big(x_n^{-\alpha}+(x_n^{-2\alpha-1+m^{-1}}-x_n^{-\alpha-1})(x-x_n)\big)\textbf{\emph{1}}
(x_n\leq x< 2x_n)\nonumber\\
&&\ \ \ \ \ \ \ \ \ \ \ \ \ \ \ \ \ \ +x_n^{-2\alpha+m^{-1}}\textbf{\emph{1}}(2x_n\leq x< x_{n+1})\Big),~x\in(-\infty,\infty).\label{201}
\end{eqnarray}
Further, let
$$\overline{G_m}(x)=(\overline{F}(x))^m=\overline{F}^m(x),~x\in(-\infty,\infty).$$
Then $G_m\in({\cal{S}}\cap{\cal{OS^*}})\setminus{\cal{S^*}}$.
\end{exam}

\begin{exam}\label{e202}
Let $m\geq 1$ be any integer. Choose any constant $\alpha\in
(2+2m^{-1},\infty)$ and any constant $x_1>4^{\alpha}$. And, for all
integers $n\geq1$, let $x_{n+1}=x_n^{1+{\alpha}^{-1}}$. Clearly,
$x_{n+1}>4x_n$ and $x_n\to\infty$ as $n\to\infty$. Now, define the
distribution $F$ as follows:
\begin{eqnarray}
\overline{F}(x)&=&\textbf{\emph{1}}(x< 0)+(x_1^{-1}(x_1^{-\alpha}-1)x^{2^{-1}}+1)\textbf{\emph{1}}(0\leq x< x_1^2)\nonumber\\
&&\ \ \ \ \ \ \ +\sum\limits_{n=1}^{\infty}
\Big(\big(x_n^{-\alpha}+(x_n^{-\alpha-2}-x_n^{-\alpha-1})(x^{2^{-1}}-x_n)\big)
\textbf{\emph{1}}(x_n^2\leq x< 4x_n^2)\nonumber\\
&&\ \ \ \ \ \ \ \ \ \ \ \ \ \ +x_n^{-\alpha-1}\textbf{\emph{1}}(4x_n^2\leq x< x_{n+1}^2)\Big),~x\in(-\infty,\infty).\label{202}
\end{eqnarray}
Further, let $G_m$ be the same as in  Example \ref{e201}, then
$G_m\in({\cal{L}}\cap{\cal{OS}^*})\setminus{\cal{S}}$, thus $G_m\notin
{\cal{S}^*}$.
\end{exam}

Therefore, the proposition is proved.\hfill$\Box$

\subsection{ The relation between the classes ${\cal{L}}_{\Delta_T}\cap{\cal{OS}}_{\Delta_T}$ and ${\cal{S}}_{\Delta_T}$}

\begin{pron}\label{p202}
For all $0<T\leq \infty$, the inclusion relation
${\cal{S}}_{\Delta_T}\subset{\cal{L}}_{\Delta_T}\cap{\cal{OS}}_{\Delta_T}\subset{\cal{L}}_{\Delta_T}\cap{\cal{OS}}$
is proper.
\end{pron}

\proof~When $T=\infty,$ the corresponding counterexamples showing
that there exist some distributions belonging to the class
${\cal{L}}\cap{\cal{OS}}\setminus{\cal{S}}$ may be found in Leslie
(1989), Lin and Wang (2012), Example \ref{e202} and Example \ref{e205} below.
So we only prove the result in the case
that $0<T<\infty.$ First, we prove a simple fact that
${\cal{L}}_{\Delta_T}\cap{\cal{OS}}_{\Delta_T}\subset{\cal{L}}_{\Delta_T}\cap{\cal{OS}}$.
Let $V\in {\cal{L}}_{\Delta_T}\cap{\cal{OS}}_{\Delta_T}$ for some
$0<T<\infty$, we have

\begin{eqnarray*}
\int_0^x\overline{V}(x-y)dV(y)&\le&\sum\limits_{k=0}^{[xT^{-1}]}
\sum\limits_{l=0}^{\infty}V(x+lT-(k+1)T+\Delta_T)V(kT+\Delta_T)\nonumber\\
&=&O\Big(\sum\limits_{l=0}^{\infty}\sum\limits_{k=0}^{[xT^{-1}]}
\int_{kT}^{(k+1)T}V(x+lT-y+\Delta_T)dV(y)\Big)\nonumber\\
&=&O\Big(\sum\limits_{l=0}^{\infty}\int_0^xV(x+lT-y+\Delta_T)dV(y)\Big)=O(\overline{V}(x)),
\end{eqnarray*}
thus $V\in\cal{OS}$, where $a(x)=O(b(x))$ mean that $\limsup a(x)(b(x))^{-1}\le1$ for two positive functions $a$ and $b$.

The following Example \ref{e205} shows that the distribution class
${\cal{L}}_{\Delta_T}\cap{\cal{OS}}_{\Delta_T}$ is properly included in the distribution class ${\cal{L}}_{\Delta_T}\cap \cal{OS}$. Now, we give three counterexamples to show that the
inclusion relationship ${\cal{S}}_{\Delta_T}\subset{\cal{L}}_{\Delta_T}\cap{\cal{OS}}_{\Delta_T}$
is also proper. To this end, we first introduce
a concept of distribution. For some distribution $F$ with a finite
and positive mean $EX$, we say that the distribution $F^{I}$ defined by
$$F^I(x)=(EX)^{-1}\int_{0}^{x}\overline{F}(y)dy\emph{\textbf{\emph{1}}}(x>0), ~~x\in(-\infty,\infty)$$
is the integrated tail distribution (or equilibrium distribution) of the distribution $F$. Related work on the integrated tail distribution can be found in Kl\"{u}ppelberg (1988), Korshunov (1997), Li and Xu (2008), and so on.

\begin{exam}\label{e203}
For any $m\geq 1$, let $G_m$ be the same as in Example \ref{e201} {or Example
\ref{e202}.} Then for all $0<T<\infty,$
$G_m^I\in({\cal{L}}_{\Delta_T}\cap{\cal{OS}}_{\Delta_T})\setminus{\cal{S}}_{\Delta_T}.$
\end{exam}

In order to give the third counterexample, we first introduce some
relevant notions and notations. We say that a distribution $F$
belongs to the exponential distribution class with the index
$\gamma\geq 0$, denoted by $F\in{\cal L}(\gamma)$, if for all
$t\in(-\infty,\infty)$, $$\overline {F}(x+t)\sim e^{-\gamma
t}\overline {F}(x).$$ We say that a distribution $F$ belongs to the
convolution equivalent distribution class with the index $\gamma\geq
0$, denoted by $F\in{\cal S}(\gamma)$, if $F\in{\cal L}(\gamma)$,
$M_{\gamma}(F)=\int_0^{\infty}e^{\gamma y}dF(y)<\infty$ and
$$\overline{F^{*2}}(x)\sim 2M_{\gamma}(F){\overline F}(x).$$
Obviously, when $\gamma=0$, ${\cal{L}}(0)$=$\cal{L}$ and
${\cal{S}}(0)$=$\cal{S}$; when $\gamma>0$, the distributions in
${\cal L}(\gamma)$ are light-tailed. The classes $\mathcal{L}(\gamma)$ and $\mathcal{S}(\gamma)$ were introduced by Chover et al. (1973, a, b) for $\gamma>0$. Bertoin and Doney (1996) note that, however, in definitions of the class $\mathcal{L}(\gamma)$ and the class $\mathcal{S}(\gamma)$, if $\gamma>0$, and if the distribution $F$ is lattice, then $x$ and $T$ should be restricted to values of the lattice span.

Further, for a distribution $F$, if $M_{\gamma}(F)<\infty$ for
some $\gamma>0$, we may define a new distribution as follows.
\begin{eqnarray*}
F_{\gamma}(x)=(M_{\gamma}(F))^{-1}\int_0^{x}e^{\gamma
y}dF(y)\textbf{1}(x\ge0),~x\in(-\infty,\infty),
\end{eqnarray*}
which is called the $\gamma$-transform or the Escher transform of
the distribution $F$. {Similarly, we can define the
$-\gamma$-transform of a distribution $F$ for any $\gamma>0$.}

\begin{exam}\label{e204}
Kl\"{u}ppelberg and Villasenor (1991) { found} two distributions
$F_i\in{\cal S}(\gamma)$ for some $\gamma>0,~i=1,2,$ but
$F=F_1*F_2\in{\cal L}(\gamma)\setminus{\cal S}(\gamma)$. Then for
all $0<T<\infty$, $F_{\gamma}\in({\cal L}_{\Delta_T}\cap{\cal
OS}_{\Delta_T})\setminus{\cal S}_{\Delta_T}$.\hfill$\Box$
\end{exam}

\begin{re}
We point out that the distribution $F_{\gamma}\in {\cal
L}_{\Delta_T}\setminus{\cal S}_{\Delta_T}$ for all $0<T<\infty$ in
Example \ref{e204} was firstly introduced by Proposition 2.1 of Chen
et al. (2013). In addition, there Example \ref{e203} and Example \ref{e204}
give two new ways to find more distributions in the class
$(\cal{L}\cap\cal{OS})\setminus{\cal S}$ and its subclasses.
\end{re}

\subsection{ The relation between the classes  $\cal{L}\cap{\cal{OS}}$ and
$\cal{L}\cap{\cal{OS}^*}$.}

\begin{pron}
The inclusion relation $
\cal{L}\cap{\cal{OS^*}}\subset\cal{L}\cap{\cal{OS}}$ is proper.
\end{pron}
\proof First, using the method of Lemma 9 of Denisov et al. (2004), we can prove the fact that $\cal{L}\cap{\cal{OS^*}}\subset\cal{L}\cap{\cal{OS}}$.
Next, we prove that the above inclusion relation is proper by
using the following example. To this end, we recall a distribution
in the distribution class $ \cal{L}\cap{\cal{OS}}\backslash\cal{S}$,
which was found by Lin and Wang (2012).
\begin{exam}\label{e205}
Let $x_1>1$ be any given number, and let $x_{n+1}=(2x_n)^2$,
$n\geq1$. For any $\alpha\in (0,1)$, define
\begin{eqnarray}
\overline{F}(x)&=&\textbf{\emph{1}}(x< 0)+(x_1^{-1}(x_1^{-\alpha}-1)x+1)\textbf{\emph{1}}(0\leq x< x_1)\nonumber\\& &\ \ \ \ \ \ \ \ \ +\sum\limits_{n=1}^{\infty}\Big(\big(x_n^{-\alpha}+(2^{-2\alpha}x_n^{-2\alpha-1}-x_n^{-\alpha-1})(x-x_n)\big)
\textbf{\emph{1}}(x_n\leq x< 2x_n)\nonumber\\
& &\ \ \ \ \ \ \ \ \ \ \ \ \ \ \ \ \ \ +(2x_n)^{-2\alpha}\textbf{\emph{1}}(2x_n\leq x< x_{n+1})\Big),~x\in(-\infty,\infty).\label{227}
\end{eqnarray}
For any positive integer $m\in(\alpha^{-1},2\alpha^{-1})$, let
$\overline{G_m}(x)=\overline{F}^m(x),~~x\in(-\infty,\infty).$ It is
obvious that the distribution $G_m$ has a finite mean. Lin and Wang
(2012) has proved that
$G_m\in(\cal{L}\cap{\cal{OS}})\backslash\cal{S}$. Further, we
have $G_m\notin\cal{OS}^*$, $G_m^I\in(\cal{L}\cap{\cal{OS}})\setminus{\cal{S}}$ and
$G_m^I\in{\cal{L}}_{\Delta_T}\setminus{\cal{OS}}_{\Delta_T}$ for all $0<T<\infty$.\hfill$\Box$
\end{exam}

\section{Local asymptotic estimations} \setcounter{equation}{0}
\setcounter{thm}{0}\setcounter{Corol}{0}\setcounter{lemma}{0}\setcounter{re}{0}

In this section, we try to deliver local asymptotic estimations
for the supremum of a random walk, where the distributions
of the summands of the random walk belong to the class
$\cal L\cap{\cal{OS^*}}$, or equivalently, the integrated tail
distributions of the summands belong to the class ${\cal
L}_{\Delta_T}\cap{\cal{OS}}_{\Delta_T}$ for any $0<T<\infty$, {see
Lemma \ref{l35} below}. This explains that distributions from the
classes $\cal L\cap{\cal{OS^*}}$ {and  ${\cal
L}_{\Delta_T}\cap{\cal{OS}}_{\Delta_T}$} possess good properties,
thus they have important value in applications.

In the following, we first introduce some concepts of a random walk
and the main result of this paper. Then  we give some lemmas in
the second subsection. The proof of the main result will be
presented at last.

\subsection{Related concepts and main result}

Let $\{X_i: i\geq 1\}$ be a sequence of independent, identically
distributed r.v.s with a common non-degenerate distribution $F$
supported on $(-\infty,\infty)$. Denote the random walk by
$\{S_n=\sum_{i=1}^{n}X_i:n\ge0\}$, where $S_0=0$,  and the supremum
of the random walk by $M=\sup_{n\ge0}S_n$ with a distribution $W$
{supported on $[0,\infty)$}. Assume that $$-\infty
<-\mu=EX_1<0,$$ then we know that $S_n$ drifts to $-\infty$ and $W$ is a proper
distribution.

Further, let $\tau_+=\inf\{n\ge1:S_n>0\}$ be the first ascending
ladder-epoch and $S_{\tau_+}$ the first ascending ladder hight with
a defective distribution $F_+$ {supported on $[0,\infty)$}, i.e.,
$0<p=F_+(\infty)<1$. Denote $G(x)=p^{-1}F_+(x),
x\in(-\infty,\infty)$, then $G$ is a proper distribution {supported
on $[0,\infty)$}. It is well known that, for any $0<T\le\infty$ and
$x\ge0$,
\be\label{31}W(x+\d)=(1-p)\sum_{n=1}^{\infty}p^nG^{*n}(x+\d),\ee
see Asmussen et al. (2002) or Asmussen et al. (2003).
\medskip

For the random walk $\{S_n:n\ge0\}$, if $F\in\cal S^*$, then
\be\label{32}W(x+\d)\sim\mu^{-1}T\overline F(x),\ee
for any $0<T<\infty$, see Asmussen et al. (2002), and so on. Further, Asmussen et al. (2003) show that $W\in{\cal S}_{\Delta_T}$ for any $0<T<\infty$. Naturally, one hopes to know that if $F\in\cal L\cap{\cal{OS^*}}$,
then how to estimate $W(x+\d)$? And, what is distribution of the supremum $M$? Our answer is as follows.

\begin{thm}\label{t31}
For the random walk $\{S_n:n\ge0\}$, if $F\in\cal L$,
then for any $0<T<\infty$,
\be\label{33}\liminf W(x+\d)(\overline F(x))^{-1}=\mu^{-1}T.\ee
Further, if $F\in\cal
L\cap{\cal{OS^*}}$ and
\be\label{34}C^\otimes(F)<\mu+2EX_{1}^{+},\ee
then for any $0<T<\infty$,
\be\label{35}
\limsup W(x+\d)\big(\overline F(x)\big)^{-1}\le\mu^{-1}T\lo(1-\mu^{-1}(C^\otimes(F)-2EX_{1}^{+})\ro)^{-1},
\ee
and $W\in{\cal L}_{\Delta_T}\cap{\cal{OS}}_{\Delta_T}$. In
particular, if $C^\otimes(F)=2EX_{1}^{+}$, namely $F\in{\cal S}^*$,
then (\ref{32}) holds. In addition, if $F\in(\cal L\cap{\cal{OS^*}})\setminus \cal {S}^*$,
then $W\in({\cal L}_{\Delta_T}\cap{\cal{OS}}_{\Delta_T})\backslash{\cal{S}}_{\Delta_T}$.
\end{thm}

\begin{re}
Here we note that the theorem gives also us a new way to find more distributions in
the class $({\cal L}_{\Delta_T}\cap{\cal{OS}}_{\Delta_T})\backslash
{\cal{S}}_{\Delta_T}$ for any $0<T<\infty$.

In addition, 
we give the following results to illustrate the condition (\ref{34})
\end{re}

\begin{pron}\label{p31}
There exists a long-tailed and generalized strong subexponential distribution $F$ supported on $(-\infty,\infty)$ satisfying the condition (\ref{34}).
\end{pron}

\subsection{Some lemmas}
In this section, we prepare more lemmas on
the local distributions, which will be used in the proof of Theorem
3.1 and also have their own independent value.

First, we recall a known fact. If a distribution
$V\in{\cal{L}}_{\Delta_T}$ for some $0<T\le \infty$, then
\begin{eqnarray*}
\mathcal{H}_{\Delta_T}(V)&=&\{h ~{\rm{on}}~
[0,\infty):~h(x)\uparrow\infty,~h(x)=o(x)\ \mbox{and}\\
&&V(x+t+\Delta_T)\sim V(x+\Delta_T)\ \mbox{holds uniformly for
all}~|t|\leq h(x)\}\neq\emptyset.
\end{eqnarray*}
And if $h\in\mathcal{H}_{\Delta_T}(V)$ and $h(x)\ge
h_1(x)\uparrow\infty$, then $h_1\in\mathcal{H}_{\Delta_T}(V)$ too.
Particularly, when $T=\infty$, we denote
$\mathcal{H}_{\Delta_{\infty}}(V)$ by $\mathcal{H}(V)$.

\begin{lemma}\label{l31}(i) If $V\in{\cal L}_{\d}$ for some $0<T\le \infty$, and for some $h\in\mathcal{H}_{\Delta_T}(V)$,
\be\label{36}\int_{h(x)}^{x-h(x)}V(x-y+\d)dV(y)=O(
V(x+\Delta_T)),\ee then  $V\in{\cal OS}_{\d}$.

(ii) If $V\in{\cal{L}}_{\Delta_T}\cap{\cal OS}_{\d}$, then for all
$h\in\mathcal{H}_{\Delta_T}(V)$, {(\ref{36}) holds and}
\be\label{37}\limsup\int_{h(x)}^{x-h(x)}V(x-y+\d)(V(x+\Delta_T))^{-1}dV(y)=C^T(V)-2.\ee\end{lemma}

\proof  By $V\in{\cal L}_{\d}$ and a standard method, we have
\be\label{38}V^{*2}(x+\d)\sim
2V(x+\Delta_T)+\int_{h(x)}^{x-h(x)}V(x-y+\d)dV(y),\ee thus
$V\in{\cal OS}_{\d}$ follows immediately from   (\ref{38}) and
(\ref{36}).

On the other hand, if $V\in{\cal OS}_{\d}$, {then (\ref{36}) and}
(\ref{37}) follow directly from (\ref{38}).$\hfill\Box$\medskip

\begin{lemma}\label{l32}If $V\in{\cal{L}}_{\Delta_T}\cap{\cal OS}_{\d}$ for some
$0<T\le\infty$, then for all $n\ge1$,
\be\label{39}\limsup V^{*n}(x+\d)(V(x+\Delta_T))^{-1}\le\sum_{k=0}^{n-1}(C^T(V)-1)^{n-1-k}.\ee
\end{lemma}
\proof Apparently, (\ref{39}) holds for $n=1,2$. We assume that
(\ref{39}) holds for $n=m$ and aim to show that it holds for $n=m+1$
too. 
For any $n\ge 1$, denote
$$C_n^T(V)=\limsup V^{*n}(x+\d)(V(x+\Delta_T))^{-1}.$$ For any
$h\in\mathcal{H}_{\Delta_T}(V)\cap \mathcal{H}_{\Delta_T}(V^{*m})$,
by a standard method, we obtain \be\label{q10}V^{*(m+1)}(x+\d)\sim
V(x+\Delta_T)+V^{*m}(x+\d)+I(x),\ee where
\bea\label{q11}I(x)&\le&\int_{h(x)-T}^{x-h(x)+T}V^{*m}(x-y+\d)dV(y)\nonumber\\&\lesssim&
C_m^T(V)\int_{h(x)-T }^{x-h(x)+T}V(x-y+\d)dV(y).\eea It follows from
(\ref{39})-(\ref{q11}) and Lemma \ref{l31} that \beq
C_{m+1}^T(V)&\le&
C_m^T(V)(C^T(V)-1)+1\\&\le&\sum_{k=0}^{m}(C^T(V)-1)^{m-k}, \eeq
namely (\ref{39}) holds for $n=m+1$. $\hfill\Box$\medskip

\begin{lemma}\label{lem5}
~~Let $V$ be a proper distribution on $[0,\infty)$. If
$V\in{\cal{L}}_{\Delta_T}\cap{\cal{OS}}_{\Delta_T}$ for some
$0<T\le\infty$, then for arbitrary $\varepsilon>0 $, there exists
$x_1>0$ and $K=K(\varepsilon, x_1)>0$ such that for all $x\geq x_1$
and $n\geq1$,
\begin{eqnarray}
{V^{*n}}(x+\Delta_T)\leq K(C^T(V)-1+\varepsilon)^n
{V}(x+\Delta_T).\label{tj1}
\end{eqnarray}
\end{lemma}

\proof For any $h\in\mathcal{H}_{\Delta_T}(V)$, by the condition
$V\in{\cal{L}}_{\Delta_T}\cap{\cal{OS}}_{\Delta_T}$, we know that
for any $\varepsilon\in (0,1)$, there exists a positive number $x_1$
large enough such that when $x\ge x_1$, $h(x)>2T$,

\begin{equation}
\left|{V}(x-y+\Delta_T)({V}(x+\Delta_T))^{-1}-1\right|<8^{-1}\varepsilon
\quad{\rm{uniformly ~for}}\quad|y|\leq h(x)\label{305}
\end{equation}
and
\begin{equation}
\int_{h(x)-T}^{x-h(x)+T}{V}(x-y+\Delta_T)dV(y)<\lo(C^T(V)-
2+8^{-1}\varepsilon\ro) {V}(x+\Delta_T).\label{306}
\end{equation}
We now prove the lemma by induction. When $n=1$, (\ref{tj1}) is
obvious. Assume that (\ref{tj1}) holds for a fixed integer $n\geq1$,
we now show that it holds for $n+1$. Denote
\begin{equation*}
A_n=\sup_{x\geq
x_1}{V^{*n}}(x+\Delta_T)({V}(x+\Delta_T))^{-1},~n\geq 1.
\end{equation*}
We have
\begin{eqnarray}
V^{*(n+1)}(x+\Delta_T)&=&\int_0^{h(x)}{V^{*n}}(x-y+\Delta_T)dV(y)
+\int_{0}^{h(x)}V(x-y+\Delta_T)dV^{*n}(y)\nonumber\\
& &+P(S_{n+1}\in x+\Delta_T, S_n>h(x), X_{n+1}>h(x))\nonumber\\
&=&I_1(x)+I_2(x)+I_3(x).\label{308}
\end{eqnarray}
{When $x\ge x_1$, by (\ref{305}), we have}
\begin{eqnarray}
I_1(x)({V}(x+\Delta_T))^{-1}&=&\int_{0}^{h(x)}
{{V^{*n}}(x-y+\Delta_T)}({V}(x+\Delta_T))^{-1}dV(y)\nonumber\\
&\leq& A_n\int_{0}^{h(x)}
{{V}(x-y+\Delta_T)}({V}(x+\Delta_T))^{-1}dV(y)\nonumber\\
&\leq& A_n\lo(1+8^{-1}\varepsilon\ro).\label{309}
\end{eqnarray}
Similarly, we have
\begin{eqnarray}
I_2(x)({V}(x+\Delta_T))^{-1}&\leq&1+8^{-1}\varepsilon.\label{310}
\end{eqnarray}
Finally, by (\ref{306}),
\begin{eqnarray}
I_3(x)({V}(x+\Delta_T))^{-1}
&\leq&\int_{h(x)-T}^{x-h(x)+T}{V^{*n}}(x-y+\Delta_T)({V}(x+\Delta_T))^{-1}dV(y)\nonumber\\
&\leq&A_n\int_{h(x)-T}^{x-h(x)+T}V(x-y+\Delta_T)({V}(x+\Delta_T))^{-1}dV(y) \nonumber\\
&\leq&A_n\lo(C^T(V)- 2+8^{-1}\varepsilon\ro).\label{311}
\end{eqnarray}
So when $x\geq x_1$, it follows from (\ref{308})-(\ref{311}) that
\begin{eqnarray}
A_{n+1}&\leq&\sup\limits_{x\geq
x_1}V^{*n}(x+\Delta_T)({V}(x+\Delta_T))^{-1}\nonumber\\
&\leq&1+8^{-1}\varepsilon+A_n\lo(C^T(V)-1+4^{-1}\varepsilon\ro).\label{tc1}
\end{eqnarray}
Taking $K=K(\epsilon)=:\frac{8}{3\varepsilon}$ and using
(\ref{tc1}), we get
\begin{eqnarray*}
A_{n+1}&\leq&2+
K\left(C^T(V)-1+3\cdot4^{-1}\varepsilon\right)\textcolor[rgb]{1.00,0.00,0.00}{^{n+1}}\nonumber\\
&\leq& K\left(C^T(V)-1+\varepsilon\right)^{n+1}.
\end{eqnarray*}This completes the proof of Lemma \ref{lem5}.
~~$\hfill\Box$\medskip

For a r.v. $\xi$ with a distribution $V$
and a positive and finite mean $E\xi$, denote
$$V_1(x)=\min\{1, V^I(x)E\xi\},~x\in(-\infty,\infty),$$
{which may be a defective distribution.

\begin{lemma}\label{l34} Let two distributions $V$ and $V_{1}$ be described as above.
If $V\in{\cal L}$, then
$V_1\in{\cal L}_{\d}$ for all $0<T<\infty$; on the contrary, if
$V_1\in{\cal L}_{\d}$ for some $0<T<\infty$, then $V\in{\cal L}$.
And both are able to derive the following asymptotic equivalence
formula: \be\label{q121}V_1(x+{\d})\sim \overline{V}(x)T.\ee
\end{lemma}
\proof Clearly, since $V\in{\cal L}$, $V_1\in{\cal L}_{\d}$ for all $0<T<\infty$ and
(\ref{q121}) holds. On the contrary, from the following fact that for any $y>0$ and
$x$ large enough,
\beq \overline{V}(x-y)T\le V_1(x-y-T+{\d})\sim
V_1(x+{\d})\le\overline{V}(x)T, \eeq
we know that $V\in{\cal L}$ and (\ref{q121}) holds.$\hfill\Box$\medskip

In the following, we give some new versions of Pitman's theorem.
When $T=\infty$, the result is due to Pitman (1980). To this end, we
first recall a known fact. Yu and Wang (2014) show that, if
a distribution $V\in{\cal L}\cap{\cal OS}$, then for any
$h\in\mathcal{H}(V)$,\beq
\limsup\int_{h(x)}^{x-h(x)}\overline{V}(x-y)(\overline{V}(x))^{-1}dV(y)=C^*(V)-2.\eeq
And by Lemma \ref{l31} (ii), if $V\in{\cal L}_{\d}\cap{\cal OS}_{\d}$
for some $0<T\le\infty$, then (\ref{37}) holds.
Similar to the proof of (\ref{37}), we can prove that, if a r.v. $\xi$ with distribution $V\in{\cal L}\cap{\cal {OS}^*}$, then for any
$h\in\mathcal{H}(V)$,
\be\label{q03}
\limsup\int_{h(x)}^{x-h(x)}\overline{V}(x-y)\overline{V}(y)(\overline{V}(x))^{-1}dy
=C^{\otimes}(V)-2E\xi.
\ee

\begin{lemma}\label{l35} Let $V$ and $U$ be two distributions,
then the following assertions hold.

(i) If $V\in{\cal L}_{\d}\cap{\cal
OS}_{\d}$ for some $0<T\le\infty$ and there exist two constants
$c_1$ and $c_2$ such that \be\label{q10} 0<c_1=\liminf\frac
{U(x+\d)}{V(x+\d)}\le \limsup\frac {U(x+\d)}{V(x+\d)}=c_2<\infty,\ee
then $U\in{\cal OS}_{\d}$ and \be\label{q01}C^T(U)-2c_1^{-1}c_2\le
c_1^{-1}c_2^2(C^T(V)-2).\ee Particularly, if $c_1=c_2=c_0$, then
$U\in{\cal L}_{\d}\cap{\cal OS}_{\d}$ and
\be\label{q12}C^T(U)-2=c_0(C^T(V)-2).\ee

(ii) If $V\in{\cal L\cap\cal OS^*}$ and there exist two constants
$c_1$ and $c_2$ such that \be\label{q130} -\infty<c_1=\liminf U(x+\d)(\overline{V}(x))^{-1}\le \limsup U(x+\d)(\overline{V}(x))^{-1}=c_2<\infty\ee
for some $0<T<\infty$, then $U\in{\cal OS}_{\d}$ and
\be\label{q133}c_1^2(c_2T)^{-1}(C_{\otimes}(V)-2E\xi)\le
C_T(U)-2\le C^T(U)-2\le c_2^2(c_1T)^{-1}(C^{\otimes}(V)-2E\xi).\ee
Particularly,
if $c_1=c_2=c_0$, then both are able to derive the following
equations:
\be\label{q130}C^T(U)-2=c_0T^{-1}(C^\otimes(V)-2E\xi)
~and~C_T(U)-2=c_0T^{-1}(C_\otimes(V)-2E\xi).\ee
Thus, $V\in{\cal L}\cap{\cal OS}^*$ if and only if $U\in{{\cal
L}_{\Delta_T}\cap\cal OS}_{\d}$ for all $0<T<\infty$.
Further, we have
\be\label{q13}C^T(V_1)-2=C^\otimes(V)-2E\xi~and~C_T(V_1)-2=c_0T^{-1}(C_\otimes(V)-2E\xi).\ee
\end{lemma}
\proof (i) In (\ref{308}), we take  $n=1$ and $V=U$, then for any
$h\in\mathcal{H}_{\d}(V)$,  \be\label{q1331}\limsup
I_i(x)(U(x+\d))^{-1}\le c{_1}^{-1}c_2,~i=1,2.\ee
Let r.v.s $\xi_1$ and
$\xi_2$ have distributions $V$ and $U$ respectively. Then by (\ref{q10}), we have
\beq &&\limsup\frac {I_3(x)}{U(x+\d)}\le
c_2\limsup \int_{h(x)-T}^{x-h(x)+T}\frac {V(x-y+\d)}{U(x+\d)}dU(y)\\
&=&c_2\limsup P(\xi_1+\xi_2\in x+\d, ~h(x)-T<\xi_2\le
x-h(x)+T)(U(x+\d))^{-1}\\
&\le&c_2\limsup P(\xi_1+\xi_2\in x+\d,~h(x)-2T<\xi_1\le x-h(x)+2T)(U(x+\d))^{-1}\\
&=&c_2\limsup\int_{h(x)-2T}^{x-h(x)+2T}U(x-y+\d)(U(x+\d))^{-1}dV(y)\\
&\le&c_1^{-1}c_2^2\limsup\int_{h(x)-2T}^{x-h(x)+2T}V(x-y+\d)(V(x+\d))^{-1}dV(y)\\
&=&c_1^{-1}c_2^2(C^T(V)-2).\eeq Thus, by (\ref{38}) of Lemma
\ref{l31}, (\ref{q1331}) and (\ref{37}), we know that $U\in{\cal
OS}_{\d}$ and (\ref{q01}) holds.

Particularly, if $c_1=c_2=c_0$, then
$$C^T(V)-2\le c_0^{-1}(C^T(U)-2).$$ Combining the above
inequalities and (\ref{q01}) yields the equality (\ref{q12}).

(ii) When $V\in{\cal L\cap\cal OS^*}$, similarly to (i), for some
$h\in\mathcal{H}_{\Delta_T}(U)$, we have
\be\limsup I_i(x)(U(x+\d))^{-1}\le c{_1}^{-1}c_2,~i=1,2.\ee
For $I_3(x)$, without loss of generality, we may assume that $l_1(x)=(x-2h_1(x))T^{-1}$
is {an} integer for $x$ large enough, where $h_1(x)=h(x)-T$. Thus, by (\ref{q10}), we have
\beq &&\limsup I_3(x)(U(x+\d))^{-1}\le
c_2\limsup \int_{h_1(x)}^{x-h_1(x)}\o V(x-y)(U(x+\d))^{-1}dU(y)\\
&=&c_2\limsup \sum_{k=1}^{l_1(x)}\int_{h_1(x)+(k-1)T}^{h_1(x)+kT}\overline{V}(x-y)(U(x+\d))^{-1} dU(y)\\
&\le&c_1^{-1}c_2^2\limsup\sum_{k=1}^{l_1(x)}\overline{V}(x-h_1(x)-kT)
\overline{V}(h_1(x)+(k-1)T)(\overline{V}(x))^{-1}\\
&=&c_1^{-1}c_2^2T^{-1}\limsup\int_{h_1(x)}^{x-h_1(x)}\overline{V}(x-y)\overline{V}(y)(\overline{V}(x))^{-1}dy\\
&=&c_1^{-1}c_2^2T^{-1}(C^\otimes(V)-2E\xi).\eeq Thus, $U\in{\cal
OS}_{\d}$.

On the other hand, we have
\be\label{q1332}\liminf I_i(x)(U(x+\d))^{-1}\ge c{_1}c_2^{-1},~i=1,2.\ee
And without loss of generality, we may assume that $l(x)=(x-2h(x))T^{-1}$ is {an}
integer for $x$ large enough, then
\beq &&\liminf I_3(x)(U(x+\d))^{-1}\ge
c_1\liminf \int_{h(x)}^{x-h(x)}\o V(x-y)(U(x+\d))^{-1}dU(y)\\
&=&c_1\liminf \sum_{k=1}^{l(x)}\int_{h(x)+(k-1)T}^{h(x)+kT}\overline{V}(x-y)(U(x+\d))^{-1} dU(y)\\
&\ge&c_1^2c_2^{-1}\liminf\sum_{k=1}^{l(x)}\overline{V}(x-h(x)-(k-1)T)\overline{V}(h_1(x)+(k-1)T)
(\overline{V}(x))^{-1}\\
&=&c_1^2c_2^{-1}T^{-1}\liminf\int_{h(x)}^{x-h(x)}\overline{V}(x-y)\overline{V}(y)(\overline{V}(x))^{-1}dy\\
&=&c_1^2c_2^{-1}T^{-1}(C_\otimes-2E\xi).\eeq
By the above four inequalities, we know that (\ref{q133}) holds.

From these results, the final result can be obtained
immediately.$\hfill\Box$\medskip

Now we discuss the relationship between a distribution $V$ and its
$\gamma$-transformation $V_{\gamma}$.

\begin{lemma}\label{l36} For some $\gamma>0$,
$V\in{\cal L}(\gamma)\cap{\cal OS}$ if and only if $V_{\gamma}\in{\cal L}_{\d}\cap{\cal OS}_{\d}$
for any $0<T<\infty$. And both of them imply the
following asymptotic equivalence formula:
\be\label{q1111}C^{T}(V_{\gamma})-2=(M_{\gamma}(V))^{-1}\gamma(C^{\otimes}(V)-2E\xi)=(M_{\gamma}(V))^{-1}C^{\ast}(V)-2.\ee
\end{lemma}
\proof According to Proposition 2.1 of Wang and Wang (2011), we know that
$V\in{\cal L}(\gamma)$ for some $\gamma>0$ if and only if $V_{\gamma}\in{\cal L}_{\d}$
for any $0<T<\infty$. And both of them are able to derive the
following asymptotic equivalence formula:
\be\label{q1112}
V_{\gamma}(x+\d)\sim(M_{\gamma}(V))^{-1}\gamma T
e^{\gamma x}\overline{V}(x).\ee
Select $h$ and $l$ as in the proof of Lemma 3.5. By (\ref{q1112}), we know that
\beq
&&\int_{h(x)}^{x-h(x)}\frac {\o V(x-y)\o V(y)}{\o V(x)}dy
\sim\frac{M_{\gamma}(V)}{\gamma T}\sum_{k=1}^{l(x)}\int_{h(x)+(k-1)T}^{h(x)+kT}\frac {V_{\gamma}(x-y+\d)V_{\gamma}(y+\d)}{V_{\gamma}(x+\d)} dy\\
&\sim&M_{\gamma}(V)(\gamma^{-1}\sum_{k=1}^{l(x)}
V_{\gamma}(x-h(x)-kT+\d)V_{\gamma}(h(x)+kT+\d)(V_{\gamma}(x+\d))^{-1}\\
&\sim&M_{\gamma}(V)\gamma^{-1}\int_{h(x)}^{x-h(x)}V_{\gamma}(x-y+\d)(V_{\gamma}(x+\d))^{-1}dV_{\gamma}(y)\\
&=&\gamma^{-1}\int_{h(x)}^{x-h(x)}\o V(x-y)(\o V(x))^{-1} dV(y).\eeq
Thus, $V\in\cal {OS}$ if and only if $V_{\gamma}\in{\cal OS}_{\d}$
for any $0<T<\infty$, and both of them are able to derive the (\ref{q1111}) holds. $\hfill\Box$
\medskip

Finally, we introduce Corollary 2.1 of Yu and Wang (2014).

\begin{lemma}\label{l37} Let $V\in{\cal{L}}(\gamma)\cap{\cal{OS}}$ for some $\gamma\geq0$ and
$U=\sum\limits_{n=0}^{\infty}p_nV^{*n}$, where $\{p_n:~n\geq0\}$ is
a sequence of nonnegative numbers satisfying
$\sum\limits_{n=0}^{\infty}p_n=1$. Suppose that there exists some
$\varepsilon_0>0$ such that
\begin{eqnarray}
\sum_{n=0}^{\infty}p_n(C^*(V)-M_\gamma(V)+\varepsilon_0)^n<\infty,\label{iv}
\end{eqnarray}
then $U\in{\cal{L}}(\gamma)\cap{\cal{OS}}$.

\end{lemma}

\subsection{Proofs of Theorem \ref{t31} and Proposition \ref{p31}}

In this section, we prove Theorem \ref{t31} and Proposition \ref{p31}, respectively.
\medskip

\noindent{\bf Proof of Theorem \ref{t31}}. We first prove (\ref{33}). By $F\in\cal L$, Corollary 3.1 of Wang
and Wang (2006) and Lemma 3.4, we know that
\bea\label{t2}G(x+\d)&\sim&
(1-p)p^{-1}\mu^{-1}F_1(x+\d)\nonumber\\&\sim&(1-p)p^{-1}\mu^{-1}T\o
F(x), \eea hence  $G\in{\cal{L}}_{\Delta_T}$. Thus by Corollary 1 of
Asmussen et al. (2003) and Theorem 3.1 of Yu et al. (2010), we know that
\be\label{t1}\liminf W(x+\d)(G(x+\d))^{-1}=p(1-p)^{-1}.\ee
So, (\ref{33}) follows from (\ref{t2}) and (\ref{t1}).

Next, we  prove (\ref{35}). We know from (\ref{t2}),
$F\in{\cal{L}}\cap{\cal{OS}^*}$, Lemmas 3.4 and 3.5 that
$F_1\in{\cal{L}}_{\Delta_T}\cap{\cal OS}_{\d}$ and
$G\in{\cal{L}}_{\Delta_T}\cap{\cal OS}_{\d}$ for all $0<T<\infty$.
By condition (\ref{34}) and Lemma 3.5, we have
\be\label{t3}p(C^T(G)-1)=p((1-p)p^{-1}\mu^{-1}(C^\otimes(F)-2EX_1^{+})+1)<1.\ee
It follows from  (\ref{31}),  (\ref{t3}), Lemmas \ref{lem5} and
\ref{l32} and the dominated convergence  theorem that
\bea\label{t4}
&&\limsup\frac{W(x+\d)}{G(x+\d)}\le(1-p)\sum_{n=1}^{\infty}p^n
\limsup\frac{G^{*n}(x+\d)}{G(x+\d)}\nonumber\\
&\le&(1-p)\sum_{n=1}^{\infty}p^n\sum_{k=0}^{n-1}(C^T(G)-1)^{n-1-k}\nonumber\\
&=&(1-p)\sum_{n=1}^{\infty}p^n((C^T(G)-1)^{n}-1)(C^T(G)-2)^{-1},
\eea
where, if $C^T(G)=2$, then we define
$((C^T(G)-1)^{n}-1)(C^T(G)-2)^{-1}=n$ by continuity. By (\ref{t3})
and (\ref{t4}), we obtain that
\beq\limsup\frac{W(x+\d)}{G(x+\d)}&\le&(1-p)(C^T(G)-2)^{-1}\lo(\frac{p(C^T(G)-1)}
{1-p(C^T(G)-1)}-\frac{p}{1-p}\ro)\nonumber\\
&=&p(1-p(C^T(G)-1))^{-1},\eeq thus (\ref{35}) follows from
(\ref{t2}), (\ref{t3}) and (\ref{q13}).

Now, we show that $W\in{\cal
L}_{\Delta_T}\cap{\cal{OS}}_{\Delta_T}$ for any $0<T<\infty$. By
(\ref{33}), (\ref{34}), (\ref{35}) and Lemma 3.5, we immediately get
$W\in{\cal{OS}}_{\Delta_T}$. Next, we prove  $W\in{\cal
L}_{\Delta_T}$ for any $0<T<\infty$. To this end, we denote the
$-\gamma$-transform of $W$ and $G$ by $U=W_{-\gamma}$ and
$V=G_{-\gamma}$, respectively. From (\ref{31}), we know that
\be\label{3001}M_{-\gamma}(W)=(1-p)(1-p_{1})^{-1}>1,\ee thus by
(\ref{3001}), we have
\be\label{3002}M_{\gamma}(V)=(M_{-\gamma}(G))^{-1}<1\ee
and
\be\label{3003}\overline{U}(x)=(1-p_{1})\sum_{n=1}^{\infty}p_{1}^n\overline{V^{*n}}(x),\ee
where $0<p_{1}=pM_{-\gamma}(G)<1$ for some $\gamma$ large enough. By
condition (\ref{t2}), Lemma 3.5 (ii) and (\ref{34}), we have
\be\label{3004}C^{T}(G)-2=(1-p)p^{-1}\mu^{-1}(C^{\otimes}(F)-2EX_{1}^{+})<(1-p)p^{-1}.\ee
From (\ref{3004}) and Lemma 3.6, we know that
\beq p_{1}(C^{\ast}(V)-M_{\gamma}(V))<1.\eeq
Thus, by (\ref{3003}) and Lemma 3.7, $W_{-\gamma}\in\cal L(\gamma)$. According to Proposition
2.1 of Wang and Wang (2011), we know that $W\in{\cal L}_{\Delta_T}$
for any $0<T<\infty$.

Finally, if $F\in(\cal L\cap{\cal{OS^*}})\setminus \cal {S}^*$,
then by Corollary 3.2 of Wang et al. (2007), we have $W\in({\cal L}_{\Delta_T}\cap{\cal{OS}}_{\Delta_T})\backslash{\cal{S}}_{\Delta_T}$
for any $0<T<\infty$.$\hfill\Box$
\medskip

\noindent{\bf Proof of Proposition \ref{p31}.} Let $Y$ be a random variable with distribution $G\in\cal{L}\cap\cal{OS}^*$ supported on $(-\infty,\infty)$ and finite mean $EY=-\lambda<0$. If $C^\otimes(G)<\lambda+2EY^+$, we take distribution $F=G$ and $\mu=\lambda$, then distribution $F$ satisfy the condition (\ref{34}). Otherwise, if $C^\otimes(G)\ge\lambda+2EY^+$, we set random variable $X=Y-a$ with distribution $F$ and mean $\mu=-EX$ for some $a$ large enough such that $C^\otimes(G)-2EY^+<\lambda+a$. It is easy to find that for any $h\in{\cal{H}}(F)={\cal{H}}(G)$,
\beq
C^\otimes(F)-2EX^+&=&\limsup\int_{h(x)}^{x-h(x)}\overline{F}(x-y)\overline{F}(y)dy(\overline{F}(x))^{-1}\nonumber\\
&=&\limsup\int_{h(x)}^{x-h(x)}\overline{G}(x-y)\overline{G}(y)dy(\overline{G}(x))^{-1}\nonumber\\
&=&C^\otimes(G)-2EY^+<\lambda+a=\mu,
\eeq
that is condition (\ref{34}) holds. $\hfill\Box$

\section{Some applications}
\setcounter{equation}{0}\setcounter{thm}{0}\setcounter{Corol}{0}\setcounter{lemma}{0}\setcounter{re}{0}

In this section, we give some practical applications of the results, concepts and methods of this paper to the renewal risk model and $M/G/1$ queue, respectively.

\subsection{On the ruin distribution}
In the renewal risk model, the claim sizes $Y_i,i\ge1$ are independent, identically distributed r.v.s with a common
non-degenerate distribution $F_{1}$ supported on $(0,\infty)$ and finite mean $EY_1$; the claims occur at the random instants of time $0<T_1<T_2<...$ a.s, and the inter-arrival times $Z_i,i\ge1$ are also independent, identically distributed r.v.s with a common
non-degenerate distribution $F_{2}$ supported on $[0,\infty)$ and finite mean $EZ_1$; the sequences $\{Y_i:i\ge1\}$ and $\{Z_i:i\ge1\}$ are independent of each other; the initial capital and the premium income rate are denoted by $x\ge0$ and $c>0$, respectively; finally, the net profit condition $cEZ_1>EY_1$ is required.

We denote r.v.s $X_i=Y_i-cZ_i,i\ge1$, which have a common distribution $F$ supported on $(-\infty,\infty)$ and mean $EX_1$. Then the distribution $W$ of the supremum $M$ of a random walk $\{S_n=\sum_{i=0}^n X_i:n\ge0\}$ is called the ruin distribution of the renewal risk model, and its tail distribution
$$\overline{W}=(1-p)\sum_{n=1}^{\infty}p^n\overline{G^{*n}}$$
is called the ruin probability, where $0<p=P(\tau_+<\infty)<1$ and $G=pF_+$. In particular, if $Z_1$ is exponentially distributed with intensity $(EZ_1)^{-1}$, then $p=EY_1(cEZ_1)^{-1}$ and $G=F^I_{1}$. See, for instance, Veravebeke (1977) or Embrechts et al. (1997).

Clearly, (\ref{31}) holds for any $0<T\le\infty$ in this model. Further, we assume that $F\in\mathcal{L}\cap\mathcal{OS}^*$ and condition (\ref{34}) is satisfied, then by Theorem \ref{31}, we can get the local asymptotic estimation for the ruin distribution $W$. In other words, if the local ruin probability is defined by $W(x+\Delta)$ for any $0<T<\infty$, then its asymptotic estimation is obtained.


\subsection{On the stationary distribution of the virtual waiting-time}
For the sake of simplicity, we omit the detailed description of the $M/G/1$ queue. As Kl\"{u}ppelberg (1989) pointed out, the arrival rate and the service-time distribution were denoted by $\eta$ and $F_1$ respectively, where $F_1$ has a finite mean $\mu(F_1)$. If $p=\eta\mu(F_1)<1$, then the stationary distribution $W$ of the virtual waiting-time can be written as
$$ W=(1-p)\sum_{n=0}^{\infty}p^n G^{*n},$$
where $G=F_1^I$. Thus (\ref{31}) holds for any $0<T\le\infty$.

Further, we assume that $F_1\in\mathcal{L}\cap\mathcal{OS}^*$ and $C^\otimes(F_1)<(1+p^{-1})\mu(F_1)$, then for any $0<T<\infty$, $G(x+\Delta)\sim T(\mu(F_1))^{-1}\overline{F_1}(x)$.
Therefore, by the proof of Theorem \ref{31} and Lemma \ref{l35} (ii), we can also get the local asymptotic estimation for the stationary distribution $W$ of the virtual waiting-time as follows:
\beq
Tp(1-p)^{-1}\big(\mu(F_1)\big)^{-1}&=&\liminf W(x+\Delta)\big(\overline{F_1}(x)\big)^{-1}\nonumber\\
&\le&\limsup W(x+\Delta)\big(\overline{F_1}(x)\big)^{-1}\nonumber\\
&\le& Tp\Big(1-p\big(\mu(F_1)\big)^{-1}\big(C^\otimes(F_1)-\mu(F_1)\big)\Big)^{-1}(\mu(F_1))^{-1}.
\eeq
Particularly, if $F_1\in\mathcal{S}^*$, that is $C^\otimes(F_1)=2\mu(F_1)$, then we have
$$W(x+\Delta)\big(\overline{F_1}(x)\big)^{-1}\sim Tp(1-p)^{-1}\big(\mu(F_1)\big)^{-1}.$$

In addition to the above two examples, in a number of areas of applied probability, such as branching processes, infinitely divisible distribution and so on, the distribution of the research objects can be written in the form of (\ref{31}), or the form of other compound distributions, so the results of this paper can also be applied to these areas. Here, we omit the details.

\section{Proofs of Examples}
\setcounter{equation}{0}\setcounter{thm}{0}\setcounter{Corol}{0}\setcounter{lemma}{0}

\subsection{Proof of $G_m\in{\cal{S}}\cap{\cal{OS^*}}\setminus{\cal{S^*}}$ in Example \ref{e201}.}

It is not hard to see that for all $m\geq 1$ and $x\geq x_1$, one has
{\begin{eqnarray} x^{-2\alpha+m^{-1}}\leq\overline{F}(x)\leq
2^{\alpha}x^{-\alpha}.\label{203}
\end{eqnarray}

Since $\alpha\in (m^{-1}, 1+m^{-1})$, by $\alpha>1$ and (\ref{203}), we know that
the distribution $G_m$ has a finite mean $m(G_m)=\mu$. Denote \beq
f(x)=\sum\limits_{n=1}^{\infty}(x_n^{-\alpha-1}-x_n^{-2\alpha-1+m^{-1}})\textbf{1}(x_n<x<
2x_n).\eeq

We first prove that  $G_m\notin{\cal{S^*}}$. To this end, for all
$x\geq0$, denote
\begin{eqnarray}\label{402}
H(x)&=&(\overline{G_m}(x))^{-1}\int_0^{x}\overline{G_m}(x-y)\overline{G_m}(y)dy
\end{eqnarray}

By (\ref{201}), one has
\begin{eqnarray}
&&H(2x_n)= 2(\overline{F}^m(2x_n))^{-1}\int_{x_n}^{2x_n}
\overline{F}^m(2x_n-y)\overline{F}^m(y)dy\nonumber\\
&=&2\int_{0}^{x_n}\overline{F}^m(y)(1+(x_n^{\alpha-m^{-1}-1}-x_n^{-1})y)^mdy\nonumber\\
&=&2\int_{0}^{x_n}\overline{F}^m(y)dy+2\int_{0}^{x_n}\overline{F}^m(y)
\Big(\big(1+(x_n^{\alpha-m^{-1}-1}-x_n^{-1})y\big)^m-1\Big)dy.\label{206}
\end{eqnarray}
By (\ref{203}) and (\ref{201}),
\begin{eqnarray}
\lim\limits_{n\to\infty}x_n^{\alpha
m-m-1}\int_{0}^{x_n}\overline{F}^m(y)y^mdy
&=&\lim\limits_{n\to\infty}x_n^{\alpha m-m-1}\int_{2x_n-1}^{x_n}\overline{F}^m(y)y^mdy\nonumber\\
&=&(m+1)^{-1}. \label{207}
\end{eqnarray}
And for all $t=1,\cdots,m-1$,
\begin{eqnarray}
\lim\limits_{n\to\infty}x_n^{(\alpha-m^{-1}-1)t}\int_{0}^{x_n}\overline{F}^m(y)y^tdy=0.\label{208}
\end{eqnarray}
By (\ref{206})-(\ref{208}),
\begin{eqnarray}
\lim_{n\to\infty}H(2x_n)=2\mu+2(m+1)^{-1},\label{209}
\end{eqnarray}
thus $G_m\notin{\cal{S^*}}$.

Next, we prove $G_m\in{\cal{OS^*}}$.
We estimate $H(x)$ in the cases $x_n\leq x<\frac{3}{2}x_n$,
$\frac{3}{2}x_n\leq x< 2x_n$ and $2x_n\leq x< x_{n+1},~n\ge 1$,
respectively. When $x\in[x_n,3\cdot2^{-1}x_n)$, by (\ref{203}) and
(\ref{201}),
\begin{eqnarray}
H(x)&\leq&2^{m\alpha+1}\lo(\overline{F}^m\lo({3\cdot2^{-1}x_n}\ro)\ro)^{-1}
\int_{\frac{x}{2}}^{x}\overline{F}^m(x-y)y^{-m\alpha}dy\nonumber\\
&\leq&2^{2m\alpha+1+m}\int_{0}^{3\cdot4^{-1}x_n}\overline{F}^m(y)dy\le2^{2m\alpha+1+m}\mu.\label{211}
\end{eqnarray}
When $x\in[3\cdot2^{-1}x_n,2x_n)$, by (\ref{201}), (\ref{206}) and (\ref{209}),

\begin{eqnarray}
H(x)&=& 2(\overline{F}^m(x))^{-1}\Big(\int_{\frac{x}{2}}^{x_n}+\int_{x_n}^{x}\Big)
\overline{F}^m(x-y)\overline{F}^m(y)dy\nonumber\\
&=&2(\overline{F}^m(x))^{-1}\Big(x_n^{-m\alpha}\int_{x-x_n}^{2^{-1}x}\overline
F^m(x-y)\overline F^m(y)dy+\int_{0}^{x-x_n}\overline
F^m(x-y)\overline F^m(y)dy\Big)\nonumber\\
&\leq&2(\overline{F}^m(2x_n))^{-1}x_n^{-m\alpha}\int_{2^{-1}x_n}^{x_n}\overline{F}^m(y)dy\nonumber\\
& &+2\int_0^{x-x_n}\overline{F}^m(y)\Big(1+(x_n(1-x_n^{-\alpha+m^{-1}})^{-1}
-(x-x_n))^{-1}y\Big)^mdy\nonumber\\
&\leq&1+2\int_0^{x_n}\overline{F}^m(y)
\Big(1+(x_n^{\alpha-m^{-1}-1}-x_n^{-1})y\Big)^mdy\nonumber\\
&=&1+H(2x_n)\rightarrow1+2\mu+2(m+1)^{-1},~~n\rightarrow\infty.\label{212}
\end{eqnarray}
When $x\in[2x_n,x_{n+1})$, by (\ref{201}) and (\ref{206}),
\begin{eqnarray}
H(x)&\leq& 2(\overline{F}^m(x))^{-1}\Big(\int_{x_n}^{2x_n}
+\int_{2x_n}^{x}\Big)\overline{F}^m(x-y)\overline{F}^m(y)dy\nonumber\\
&\leq&H(2x_n)+2\int_0^{x-2x_n}\overline{F}^m(y)dy
\rightarrow4\mu+2(m+1)^{-1},~~~n\rightarrow\infty.\label{213}
\end{eqnarray}
It follows from (\ref{211})-(\ref{213}) that $G_m\in{\cal{OS^*}}$.

Finally, we prove $G_m\in{\cal{S}}$. By (\ref{203}), we have
\begin{eqnarray}
\overline{F}^{2m}(2^{-1}x)(\overline{F}^m(x))^{-1}
&\leq&2^{4m\alpha}x^{-1}\to0.\label{214}
\end{eqnarray}
By (\ref{214}) and
\begin{eqnarray*}
\overline{G_m^{*2}}(x)=2\overline{G_m}(x)-\overline{G_m}^2(2^{-1}x)
+2\int_{2^{-1}x}^{x}\overline{G_m}(x-y)dG_m(y),
\end{eqnarray*}
we know that in order to prove  $G_m\in{\cal{S}}$, it suffices to
prove
\begin{eqnarray}
T(x)&=&2(\overline{F}^m(x))^{-1}\int_{2^{-1}x}^{x}\overline{F}^m(x-y)d\lo(1-\overline{F}^m(x)\ro)\nonumber\\
&=&2m(\overline{F}^m(x))^{-1}\int_{2^{-1}x}^{x}\overline{F}^m(x-y)\overline{F}^{m-1}(y)f(y)dy\to0.\label{215}
\end{eqnarray}
Clearly, $T(x_n)=0,n\geq1$. By (\ref{201}) and (\ref{208}), we have
\begin{eqnarray}
&&T(2x_n)=2m\overline{F}^m(2x_n)\int_{0}^{x_n}\overline{F}^m(y)
\overline{F}^{m-1}(2x_n-y)f(2x_n-y)dy\nonumber\\
&=&2m(x_n^{\alpha-m^{-1}-1}-x_n^{-1})\int_{0}^{x_n}
\overline{F}^m(y)(1+(x_n^{\alpha-m^{-1}-1}-x_n^{-1})y)^{m-1}dy\nonumber\\
&\leq&2mx_n^{\alpha-m^{-1}-1}\int_{0}^{x_n}\overline{F}^m(y)
(1+x_n^{\alpha-m^{-1}-1}y)^{m-1}dy\nonumber\\
&=&2mx_n^{\alpha-m^{-1}-1}\Big(\int_{0}^{x_n}\overline{F}^m(y)dy+
\int_{0}^{x_n}\overline{F}^m(y)((1+x_n^{\alpha-m^{-1}-1}y)^{m-1}-1)dy\Big)\nonumber\\
&\to&0, ~~n\to\infty.\label{217}
\end{eqnarray}
In the following, we prove (\ref{215}) in the cases $x_n\leq x< 2x_n$
and $2x_n\leq x< x_{n+1},~n\ge 1$, respectively. When $x\in[x_n,2x_n)$,
by (\ref{201}) and (\ref{217}),
\begin{eqnarray}
&&T(x)=  2m(\overline{F}^m(x))^{-1}\int_{x_n}^{x}\overline{F}^m(x-y)
\overline{F}^{m-1}(y)f(y)dy\nonumber\\
&=&2m(\overline{F}^m(x))^{-1}\int_{0}^{x-x_n}\overline{F}^m(y)
\overline{F}^{m-1}(x-y)f(x-y)dy\nonumber\\
&\leq& 2m(x_n^{\alpha-m^{-1}-1}-x_n^{-1})\int_{0}^{x-x_n}
\overline{F}^m(y)\Big(1+(x_n(1-x_n^{-\alpha+m^{-1}}
)^{-1}-(x-x_n))^{-1}y\Big)^{m-1}dy\nonumber\\
&\leq& T(2x_n)\to0, ~~n\to\infty.\label{218}
\end{eqnarray}
When $x\in[2x_n,x_{n+1})$, by (\ref{201}) and (\ref{217}),
\begin{eqnarray}
T(x)&\leq&2m(\overline{F}^m(2x_n))^{-1}\int_{x_n}^{2x_n}
\overline{F}(2x_n-y)\overline{F}^{m-1}(y)f(y)dy\nonumber\\
&=& T(2x_n)\to0, ~~n\to\infty.\label{219}
\end{eqnarray}
According to (\ref{218}) and (\ref{219}), (\ref{215}) holds, thus
$G_m\in{\cal{S}}$.

In summary, we have
$G_m\in({\cal{S}}\cap{\cal{OS^*}})\setminus{\cal{S^*}}$.\hfill$\Box$
\medskip

\subsection{ Proof of
$G_m\in({\cal{L}}\cap{\cal{OS}^*})\setminus{\cal{S}}$ in Example
\ref{e202}.}

Still let $f$ be the density of $F$, when $x\geq x_1,$
it is easily seen that
\begin{eqnarray}
x^{-2^{-1}(\alpha+1)}\leq\overline{F}(x)\leq
2^{\alpha}x^{-2^{-1}\alpha}
\mbox{~and~}f(x)\leq2^{\alpha+1}x^{-2^{-1}\alpha-1}\label{220}
\end{eqnarray}

Moreover, one can easily find that the distribution $G_m$ has a
finite mean for $m\geq 1$, and we still denote it by $\mu.$

Firstly, by (\ref{220}) we have,
\begin{eqnarray*}
(\overline F(x))^{-1}f(x)\leq2^{\alpha+1}x^{-2^{-1}}\to0, ~~ ,
\end{eqnarray*}
thus $F\in{\cal{L}}$, so $G_m\in{\cal{L}}$. Next, we prove that
$G_m\in{\cal{OS}^*}.$ Let $H(x), x\geq0$ be the same as in
(\ref{402}), it is easily seen that
\begin{eqnarray*}
H(4x_n^2)&=&2(\overline{F}^m(4x_n^2))^{-1}\int_{0}^{2x_n^2}\overline{F}^m(y)\overline{F}^m(4x_n^2-y)dy\nonumber\\
&=&2\int_{0}^{2x_n^2}\overline F^m(y)(1+(1-x_n^{-1})(4x_n^2-y)^{2^{-1}}+2x_n)^{-1}y)^mdy\nonumber\\
&\leq&2\int_{0}^{2x_n^2}\overline F^m(y)(1+(2x_n)^{-1}y)^mdy\nonumber\\
&\leq&2x_1^2(1+(2x_n)^{-1}x_1^2)^m+2^{m\alpha+1}
\int_{x_1^2}^{\infty}y^{m-\frac{m\alpha}{2}}(y^{-1}+(2x_n)^{-1})^mdy\nonumber\\
&\leq&2x_1^2(1+x_1^2)^m+2^{m\alpha+1}(2^{-1}m\alpha-m-1)^{-1}x_1^{2m-m\alpha}(1+x_1^2)^m<\infty.
\end{eqnarray*}
Just as in  Example \ref{e201}, we deal with $H(x)$ in three cases
$x_n^2\leq x<2x_n^2$, $2x_n^2\leq x< 4x_n^2$ and
$4x_n^2\leq x< x_{n+1}^2,~n\ge 1$, respectively. When $x\in[x_n^2,2x_n^2)$,
just as (\ref{211}), by (\ref{202}) and variable substitution, we
have
\begin{eqnarray}
H(x)&=&2(\overline F^{m}(x))^{-1}\int_{\frac{x}{2}}^{x}\overline F^{m}(x-y)\overline F^{m}(y)dy\nonumber\\
&\leq&2(\overline F^{m}(2x_n^2))^{-1}\int_{\frac{x}{2}}^{x}\overline F^{m}(x-y)x_n^{-m\alpha}dy\nonumber\\
&\leq&2(2-2^{2^{-1}})^{-m}\int_{0}^{x_n^2}\overline F^{m}(y)dy<\infty.\label{222}
\end{eqnarray}
When $x\in[2x_n^2,4x_n^2)$, just as (\ref{212}), by (\ref{202}), we
have

\begin{eqnarray}
&&H(x)=2(\overline F^{m}(x))^{-1}\int_{0}^{2^{-1}x}\overline F^{m}(y)\overline F^{m}(x-y)dy\nonumber\\
&=&2(\overline F^{m}(x))^{-1}\int_{0}^{2^{-1}x}\overline F^{m}(y)
\Big(\overline F(x)+(x_n^{-\alpha-1}-x_n^{-\alpha-2})(x^{2^{-1}}+(x-y)^{2^{-1}})^{-1}y\Big)^mdy\nonumber\\
&=&2\int_{0}^{2^{-1}x}\overline F^{m}(y)\Big(1+(x_n^{-\alpha-1}-x_n^{-\alpha-2})\big(\overline F(x)(\sqrt{x}+\sqrt{x-y})\big)^{-1}y\Big)^mdy\nonumber\\
&\leq&2\int_{0}^{2x_n^2}\overline F^{m}(y)(1+(2^{2^{-1}}x_n)^{-1}y)^mdy\nonumber\\
&\leq&2x_1^2(1+(2^{2^{-1}}x_n)^{-1}x_1^2)^m+2^{m\alpha+1}\int_{x_1^2}^{\infty}y^{m-\frac{m\alpha}{2}}(y^{-1}
+(2^{2^{-1}}x_n)^{-1})^mdy\nonumber\\
&\leq&2x_1^2(1+x_1^2)^m+2^{m\alpha+1}(2^{-1}m\alpha-m-1)^{-1}x_1^{2m-m\alpha}(1+x_1^2)^m<\infty.\label{223}
\end{eqnarray}
When $x\in[4x_n^2,x_{n+1}^2)$, just as (\ref{213}), by (\ref{202}),
we have
\begin{eqnarray}
H(x)&\leq&2(\overline F^{m}(x))^{-1}\lo(\int_{2x_n^2}^{4x_n^2}+\int_{4x_n^2}^{x}\ro)\overline F^{m}(x-y)\overline F^{m}(y)dy\nonumber\\
&\leq&H(4x_n^2)+2\int_{0}^{x-4x_n^2}\overline F^{m}(y)dy<\infty.\label{224}
\end{eqnarray}
According to (\ref{222}) and (\ref{223}), (\ref{224}) holds, that is
$G_m\in{\cal{OS^*}}$.

Finally, we prove  $G_m\notin{\cal{S}}$. Since
\begin{eqnarray*}
\overline{G_m}^2(2x_n^2)(\overline{G_m}(4x_n^2))^{-1}&=&x_n^{m-m\alpha}
\big((2-2^{2^{-1}})+(2^{2^{-1}}-1)x_n^{-1}\big)^{2m}\to0,~~n\to\infty.
\end{eqnarray*}
Therefore, in order to prove  $G_m\notin{\cal{S}}$, we only need to prove
\begin{eqnarray}
\liminf_{ }T(x)=\liminf_{x\rightarrow
\infty}(\overline{G_m}(x))^{-1}\int_{2^{-1}x}^{x}\overline{G_m}(x-y)G_m(dy)
>0.\label{225}
\end{eqnarray}
In fact,
\begin{eqnarray*}
&&T(4x_n^2)=m(\overline{F}^m(4x_n^2))^{-1}\int_{2x_n^2}^{4x_n^2}\overline{F}^m(4x_n^2-y)\overline{F}^{m-1}(y)f(y)dy\nonumber\\
&=&2^{-1}m(1-x_n^{-1})\int_{0}^{2x_n^2}\overline{F}^m(y)
\Big(1+(1-x_n^{-1})(\sqrt{4x_n^2-y}+2x_n)^{-1}y\Big)^{m-1}y^{-2^{-1}}dy\nonumber\\
&\geq&2^{-1}m(1-x_n^{-1})\int_{x_1^2}^{2x_n^2}\overline{F}^m(y)y^{-\frac{1}{2}}dy\nonumber\\
&\geq&m(m\alpha+m-1)^{-1}(1-x_n^{-1})\big(x_1^{1-m\alpha-m}-(2x_n^2)^{-2^{-1}(m\alpha+m-1)}\big)\nonumber\\
&\to& m(m\alpha+m-1)^{-1}x_1^{1-m\alpha-m}>0,~~n\to\infty,
\end{eqnarray*}
that is, (\ref{225}) holds, that is $G_m \notin{\cal{S}}$.

In summary, we have
$G_m\in({\cal{L}}\cap{\cal{OS^*}})\setminus{\cal{S}}$.\hfill$\Box$

\subsection{ Proof of
$G_m^I\in({\cal{L}}_{\Delta_T}\cap{\cal{OS}}_{\Delta_T})\setminus{\cal{S}}_{\Delta_T}$
in Example \ref{e203}.}

This conclusion follows directly from
$G_m\in({\cal{S}}\cap{\cal{OS^*}})\setminus{\cal{S^*}}$, Lemmas 3.4,
3.5 and Lemma 4.2 of Wang et al. (2007).\hfill$\Box$

\subsection{ Proof of $G_m^I\in{\cal{OS}}_{\Delta_T}\setminus{\cal{S}}_{\Delta_T}$ for all
$0<T<\infty$ in Example \ref{e204}.}

By Proposition 2.1 of  Chen et al. (2013), we have $F_{\gamma}\in
{\cal{L}}_{\Delta_T}\backslash{\cal{S}}_{\Delta_T}$ for all $T>0$.
We now prove $F_{\gamma}\in {\cal{OS}}_{\Delta_T}$.  Since
$F_i\in{\cal{S}}(\gamma)\subset{\cal{OS}},~i=1,2$, so
$F=F_1*F_2\in{\cal{OS}}$ by Proposition 6.1 of Yu and Wang (2013).
By $F\in{\cal{OS}}$ and (2.1) of  Chen et al. (2013),  we have
\begin{eqnarray*}
\limsup F_{\gamma}^{*2}(x+\Delta_T)(F_{\gamma}(x+\Delta_T))^{-1}=
(M_{\gamma}(F))^{-1}\limsup\overline{F^{*2}}(x)(\overline{F}(x))^{-1}<\infty,
\end{eqnarray*}
thus  $F_{\gamma}\in {\cal{OS}}_{\Delta_T}$.\hfill$\Box$

\subsection{Proofs of $G_m\notin\cal{OS}^*$,
$G_m^I\in(\cal{L}\cap{\cal{OS}})\setminus{\cal{S}}$ and
$G_m^I\in{\cal{L}}_{\Delta_T}\setminus{\cal{OS}}_{\Delta_T}$ for all $0<T<\infty$ in Example
\ref{e205}.}

To show that $G_m\notin\cal{OS}^*$, we denote the
density of $F$ by $f$. Consider the following quantity
\begin{eqnarray*}
I(x)=\big(\overline{G_m}(2x_n)\big)^{-1}\int_{x_n}^{2x_n-2x_{n-1}}\overline{G_m}(2x_n-y)\overline{G_m}(y)dy.\nonumber\\
\end{eqnarray*}
Since  $x_n\leq y\leq {2x_n-2x_{n-1}}\leq 2x_n$, we have
$2x_{n-1}\leq {2x_n-y}\leq x_n,~n\ge 1$, so by (\ref{227}),
\begin{eqnarray*}
I(x)&=&2^{2m\alpha}x_n^{m\alpha}\int_{x_n}^{2x_n-2x_{n-1}} (x_n^{-\alpha}-f(x_n)(y-x_n))^m dy\nonumber\\
&=&2^{2m\alpha}(m+1)^{-1}x_n^{m\alpha}(x_n-2x_{n-1})\sum_{i=0}^{m}x_n^{-\alpha i}(x_n^{-\alpha}-f(x_n)(x_n-2x_{n-1}))^{m-i}\\
&\geq&2^{2m\alpha}(m+1)^{-1}(x_n-2x_{n-1})\rightarrow\infty,~~n\rightarrow\infty,
\end{eqnarray*}
thus  $G_m\notin{\cal{OS}}^*$.

Since $G_m\in{\cal{L}}\setminus{\cal{OS}^*}$, we immediately get
$G_m^I\in{\cal{L}}_{\Delta_T}\setminus{\cal{OS}}_{\Delta_T}$ by
Lemmas 3.4 and 3.5.

Now we show that $G_m^I\notin{\cal{S}}$. By (\ref{227}), it is
obvious that

{\begin{eqnarray} x^{-2\alpha}\leq\overline{F}(x)\leq
x^{-\alpha}.\label{228}
\end{eqnarray}}
By (\ref{228}), we have
\begin{eqnarray*}
(\overline{G_m^I}(x_n))^2(\overline{G_m^I}(2x_n))^{-1}&\geq&
\mu^{-1}\Big(\int_{x_n}^{2x_n}\overline{F}^m(y)dy\Big)^2
\Big(\int_{2x_n}^{\infty}\overline{F}^m(y)dy\Big)^{-1}\nonumber\\
&\geq&\mu^{-1}(m+1)^{-2}2^{2m\alpha-2}(1-m^{-1}\alpha^{-1})>0,
\end{eqnarray*}
thus by (2.4) of  Murphree (1989), we complete the proof.

To prove  $G_m^I\in{\cal{OS}}$, we consider

\begin{eqnarray}
&&\int_{\frac{x}{2}}^{x}\overline{G_m^I}(x-y)dG^I_m(y)=
\mu^{-1}  \int_{\frac{x}{2}}^{x}\overline{F}^m(y)\int_{x-y}^{\infty}\overline{F}^m(z)dzdy\nonumber\\
&=&\mu^{-1}\bigg(\int_{\frac{x}{2}}^{x}\overline{F}^m(y)\int_{x-y}^{\frac{x}{2}}\overline{F}^m(z)dzdy
+\Big(\int_{\frac{x}{2}}^{x}\overline{F}^m(y)dy\Big)^2
+\int_{\frac{x}{2}}^{x}\overline{F}^m(y)\int_{x}^{\infty}\overline{F}^m(z)dzdy\bigg)\nonumber\\
&=&\mu^{-1}(I_1(x)+I_2(x)+I_3(x)).\label{229}
\end{eqnarray}

We first estimate $I_1(x)$ in the cases $x_n\leq x<4x_n$ and $4x_n\leq
x< x_{n+1},~n\ge 1$, respectively.  When $x\in [x_n,4x_n)$, by
(\ref{227}) and (\ref{228}), we have
\begin{eqnarray}
&&\Big(\int_{x}^{\infty}\overline{F}^m(y)dy\Big)^{-1}I_1(x)
\leq\Big(\int_{4x_n}^{x_{n+1}}\overline{F}^m(y)dy\Big)^{-1}
\int_{2^{-1}x}^{x}x_n^{-m\alpha}\int_{x-y}^{2^{-1}x}\overline{F}^m(z)dzdy\nonumber\\
&\leq&2^{2m\alpha-2}(x_n^{2-m\alpha}-x_n^{1-m\alpha})^{-1}\int_0^{2^{-1}x}
\int_{y}^{2^{-1}x}\overline{F}^m(z)dzdy\nonumber\\
&\leq&2^{2m\alpha-1}(x_n^{2-m\alpha})^{-1}\Big(\int_{x_1}^{2x_n}\int_{y}^{2x_n}z^{-m\alpha}dzdy
+\int_{0}^{x_1}\int_{x_1}^{2x_n}z^{-m\alpha}dzdy
+\int_{0}^{x_1}\int_{y}^{x_1}dzdy\Big)\nonumber\\
&\leq&2^{m\alpha+1}(m\alpha-1)^{-1}(2-m\alpha)^{-1}+((m\alpha-1)^{-1}x_1^{2-m\alpha}+x_1^2)<\infty.\label{230}
\end{eqnarray}
Similarly,  when $x\in [4x_n,x_{n+1})$, we have
\begin{eqnarray}
&&\Big(\int_{x}^{\infty}\overline{F}^m(y)dy\Big)^{-1}I_1(x)
\leq\Big(\int_{x_{n+1}}^{2x_{n+1}}\overline{F}^m(y)dy\Big)^{-1}
\int_{2^{-1}x}^{x}(2x_n)^{-2m\alpha}\int_{x-y}^{2^{-1}x}\overline{F}^m(z)dzdy\nonumber\\
&\leq&
\frac{m+1}{2x_n}\Big(\int_{x_1}^{\frac{x_{n+1}}{2}}\int_{y}^{\frac{x_{n+1}}{2}}z^{-m\alpha}dzdy+
\int_{0}^{x_1}\int_{x_1}^{\frac{x_{n+1}}{2}}z^{-m\alpha}dzdu+
\int_{0}^{x_1}\int_{y}^{x_1}dzdy\Big)\nonumber\\
&\leq&(m+1)((m\alpha-1)^{-1}(2-m\alpha)^{-1}2^{-m\alpha}
+((m\alpha-1)^{-1}x_1^{2-m\alpha}+x_1^2)<\infty.\label{231}
\end{eqnarray}
Next, we estimate $I_2(x)$. When $x\in[x_n,x_{n+1})$, by (\ref{228}), we have
\begin{eqnarray}
&&\Big(\int_{x}^{\infty}\overline{F}^m(y)dy\Big)^{-1}\Big(\int_{2^{-1}x}^{x}\overline{F}^m(y)dy\Big)^2
\leq\Big(\int_{x_{n+1}}^{2x_{n+1}}\overline{F}^m(y)dy\Big)^{-1}
\Big(\int_{2^{-1}x_n}^{x_{n+1}}\overline{F}^m(y)dy\Big)^2\nonumber\\
&\leq&(m+1)\Big(2^{-4+2m\alpha}9+3x_n^{1-m\alpha}+(2x_n)^{2-2m\alpha}\Big)<\infty.\label{232}
\end{eqnarray}
Finally, for $I_3(x)$, it is obvious that
\begin{eqnarray}
\Big(\int_{x}^{\infty}\overline{F}^m(y)dy\Big)^{-1}I_3(x)
&\leq&\int_{2^{-1}x}^{x}\overline{F}^m(y)dy\to0.\label{233}
\end{eqnarray}
By (\ref{229})-(\ref{233}), we get $G_m^I\in{\cal{OS}}$.\hfill$\Box$

\vspace{0.2cm}

{\bf Acknowledgement.}
The authors are grateful to the two referees for his/her
careful reading and valuable comments and suggestions, which greatly improve the original version of this paper.


\begin{thebibliography}{99}
\bibitem{}D. Denisov, S. Foss and  D. Korshunov. {\it Tail asymptotics for the supremum of a random
walk when the mean is not finite}. Queueing Systems, 46 (2004), 15-33.
\end{thebibliography}
\end{document}